\documentclass[12pt]{amsart}

\usepackage{amsmath}
\usepackage{amsthm}
\usepackage{amssymb}
\usepackage{mathabx}
\usepackage{mathrsfs} 
\usepackage{mathtools} 
\usepackage[all]{xy} 
\usepackage{graphicx}
\graphicspath{{figures/}} 

\usepackage{quickmath} 

\usepackage{cancel} 
\usepackage{enumerate} 
\usepackage{hyperref} 
\usepackage[usenames,dvipsnames]{xcolor} 
\usepackage[marginpar]{todo}

\usepackage{tikz}
\usetikzlibrary{arrows.meta}
\usetikzlibrary{calc}
\usetikzlibrary{decorations.pathmorphing}

\hypersetup{colorlinks=true, linkcolor=NavyBlue, urlcolor=WildStrawberry}

\title{Gilman's Conjecture}

\author{Andy Eisenberg}
\address{Department of Mathematics, Oklahoma State University, USA}
\email{andrew.eisenberg@slu.edu}

\author{Adam Piggott}
\address{Department of Mathematics, Bucknell University, USA}
\email{adam.piggott@uq.edu.au}

\date{\today}
\thanks{This work was partially supported by a grant from the Simons Foundation (\#317466 to Adam Piggott)}


\tikzstyle{CayVert}=[circle,fill,inner sep=2pt,outer sep=3pt]
\tikzstyle{CayEdge}=[-{Computer Modern Rightarrow[length=1.5mm]}]
\tikzstyle{CayDashEdge}=[-{Computer Modern Rightarrow[length=1.5mm]},dashed]
\tikzstyle{CayWord}=[-{Computer Modern Rightarrow[length=1.5mm]}, decorate, decoration=snake]

\newtheoremstyle{RetheoremStyle}	
	{}						
	{}						
	{\itshape}						
	{}								
	{\bfseries}						
	{.}								
	{ }								
	{\thmname{#1}\thmnote{ #3}}		

\theoremstyle{RetheoremStyle}
\newtheorem{Retheorem}{Theorem}
\newcommand{\rethm}{\begin{Retheorem}}
\newcommand{\erethm}{\end{Retheorem}}

\begin{document}

\maketitle

\begin{abstract}
We prove a conjecture made by Gilman in 1984 that the groups presented by finite, monadic, confluent rewriting systems are precisely the free products of free and finite groups.
\end{abstract}

\tableofcontents

\section{Introduction}
\label{sec:Introduction}

Many algebraic structures are defined by, or at least naturally accompanied by, a finite rewriting system.  A \emph{rewriting system} is a pair $(\Sigma, T)$, where $\Sigma$ is a finite alphabet of symbols, $\Sigma^\ast$ denotes the set of all words over the alphabet $\Sigma$, and $T \subset \Sigma^* \times \Sigma^*$ is a set of rewriting rules.   Each rewriting rule $(L, R)$ specifies an allowable replacement: whenever $L$ appears as a subword, it may be replaced by $R$.  We write $U \xrightarrow{*} V$, if the word $U$ can be transformed into the word $V$ by application of a finite sequence of rewriting rules.  The reflexive and symmetric closure of $\xrightarrow{*}$ is an equivalence relation on $\Sigma^*$ whose equivalence classes form a monoid under the operation of concatenation of representatives.  Sometimes this monoid is a group.

A fundamental question of combinatorial group theory and the foundations of computer science asks which algebraic classes of groups can be characterized by the types of rewriting systems presenting groups in that class.  Having a nice rewriting system for a particular group often allows one to perform efficient computations in the group---for example, solving the word or conjugacy problems. A substantial effort, with contributions from many authors spanning a period of more than three decades 
(\cite{Cochet}, \cite{GilmanConjecture}, \cite{AvenhausMadlener}, \cite{AvenhausMadlenerOtto}, \cite{AutebertBoassonSenizergues}, \cite{Diekert}, \cite{MadlenerOttoGroups}, \cite{ParkesEtAl}, \cite{GilmanHermillerHoltRees}, \cite{Piggott}, and more), has been made in pursuit of a complete algebraic characterization of groups presented by \emph{length-reducing} rewriting systems (those in which each application of a rewriting rule shortens a word).  A summary of many results in this program can be found in \cite{MadlenerOttoGroups}; we mention a few relevant results here.

One can strengthen the requirement that $(\Sigma, T)$ is length-reducing in various ways, restricting attention to \emph{monadic}, \emph{2-monadic}, or \emph{special} rewriting systems.  (See Section \ref{subsec:properties} for precise definitions.)  It is common to consider \emph{confluent} rewriting systems, but this can be relaxed to require only that a rewriting system is \emph{confluent on $[1]$}, the equivalence class of the empty word (see, for example, \cite{GilmanHermillerHoltRees}, \cite{ParkesEtAl}).  

Cochet \cite{Cochet} proved that a group $G$ is presented by a finite, special, confluent rewriting system if and only if $G$ is the free product of finitely many cyclic groups.  Diekert \cite{Diekert} showed that every group presented by a finite, monadic, confluent rewriting system is virtually free.  If, in addition, the rewriting system is \emph{inverse-closed} (every element represented by a generator has an inverse which is represented by a generator), then Avenhaus and Madlener \cite{AvenhausMadlener} showed that $(\Sigma, T)$ must present a \emph{plain group}, that is, a free product of a finitely generated free group with finitely many finite groups.  Gilman \cite{GilmanConjecture} conjectured in 1984 that this was the case even without assuming that $(\Sigma, T)$ is inverse-closed.  Avenhaus, Madlener and Otto \cite{AvenhausMadlenerOtto} proved Gilman's conjecture in the special case that in each rewriting rule the left-hand side has length exactly two.  The second author proved Gilman's conjecture in the special case that every generator has finite order \cite{Piggott}.  Our main result resolves Gilman's conjecture in its full generality:

\rethm[\ref{thm:GilmansConjecture}]
A group $G$ is presented by a finite, monadic, confluent rewriting system $(\Sigma, T)$ if and only if $G$ is a plain group.
\erethm

We also give a new proof of Cochet's result by different
methods.  (See Theorem \ref{thm:Cochet}.)

In order to complete the program laid out in \cite{MadlenerOttoGroups}, it only remains to characterize the precise class of groups presented by finite, length-reducing, confluent rewriting systems.  This class is known to contain all plain groups and be a proper subclass of virtually free groups \cite{Diekert}.  It has been conjectured that this class is also the class of plain groups.  Our arguments make essential use of strong geometric consequences of the monadic hypothesis captured in Lemma \ref{lem:PlainGeometry} and therefore do not readily extend to the length-reducing setting.

\section{Background}
\label{sec:Background}

\subsection{Notation}

Throughout what follows, $\Sigma$ is a nonempty set, $\Sigma^*$ is the set of finite length words over $\Sigma$, and $T$ is a subset of $\Sigma^* \times \Sigma^*$.  The elements of $\Sigma$ are called \emph{letters}, and $\Sigma$ is the \emph{alphabet}.  The elements of $T$ are called \emph{rewriting rules}, and the pair $(\Sigma, T)$ is a \emph{rewriting system}.  We will typically use lowercase letters late in the Roman alphabet ($x, y, z, \dotsc$) to represent single letters in $\Sigma$, while uppercase letters late in the Roman alphabet ($U, V, W, \dotsc$) will represent words in $\Sigma^*$.  We will write 1 for the empty word.

If $(L, R)$ is a rewriting rule in $T$, we will write $U \to V$ to mean that $U$ contains $L$ as a subword, and $V$ is the result of replacing that subword with $R$.  We say that \emph{$V$ is obtained from $U$ by application of the rule $(L, R)$}.  We will write $U \xrightarrow{*} V$ to mean that $V$ may be obtained from $U$ by applying a finite sequence of rewriting rules, and we extend $\xrightarrow{*}$ by taking the reflexive and symmetric closure to get an equivalence relation, $\xleftrightarrow{*}$.  We write $[U]$ for the equivalence class of $U$. The set of equivalence classes, equipped with the rule $[U][V] = [UV]$, forms a monoid $M$ with identity element $[1]$.  We say that the rewriting system $(\Sigma, T)$ presents $M$.  We shall be interested in the special case that the monoid presented by a rewriting system is a group.  This happens exactly when each equivalence class represented by a letter $[x]$ has an inverse (which may or may not be represented by a letter).

If $(\Sigma, T)$ presents a group $G$, the equivalence classes $[U]$ may be identified with the group elements.  We will typically use lowercase letters early in the Roman alphabet ($a, b, c, \dotsc$) to represent group elements.  In an equation like $wx = yz$ or $UV = WX$, we mean equal as words in $\Sigma^*$.  In an equation of the form $a := U$, we mean that $a = [U]$.  By a slight abuse of notation, we will write 1 for the identity element of $G$ (which is the equivalence class of the empty word).

If $U = x_1x_2\dotsm x_p$, then $|U| = p$ is the length of the word.  For $a\in G$, we will write $|a|$ for the length of the shortest word $U$ such that $a := U$.

In pictures of portions of Cayley graphs, we will omit brackets, but any letters or words that appear as vertex labels should be understood to refer to group elements (since the vertices of the Cayley graph are the elements of the group, not the words of $\Sigma^*$).  Labels along edges should be understood to be letters in $\Sigma$, and snaking arrows will represent paths whose length may be greater than 1, which may be labeled by words from $\Sigma^*$.  Hopefully the distinction between letters, words, and group elements will be clear from context.

\subsection{Rewriting System Properties}
\label{subsec:properties}

Suppose that $(\Sigma, T)$ is a rewriting system.  A common use of a rewriting system is to construct algorithms which find \emph{normal forms}, that is, a preferred spelling of words within a particular equivalence class.  For example, one might hope to tell whether two words $U$ and $V$ are equivalent by finding their respective normal forms, which should be the same if $U \xleftrightarrow{*} V$.  Towards that end, the following properties of rewriting systems can help guarantee that the rewriting process proceeds unambiguously and terminates in finite time.

\defi
A rewriting system $(\Sigma, T)$ is called 
\enum
\item \emph{finite} if both $\Sigma$ and $T$ are finite;
\item \emph{confluent} if, whenever $W \xrightarrow{*} U$ and $W \xrightarrow{*} V$, there exists a word $Q$ so that $U \xrightarrow{*} Q$ and $V \xrightarrow{*} Q$;  
\item \emph{terminating}, or \emph{Noetherian},  if any rewriting sequence must terminate in a finite number of steps;
\item \emph{convergent} if it is both confluent and terminating.  
\eenum
\edefi

A word $U$ to which no rewriting rule can be applied is called \emph{reduced} or \emph{irreducible}, and it is clear from the definitions that an equivalence class of words in a convergent rewriting system $(\Sigma, T)$ contains a unique irreducible word.  Moreover, given any word, we may apply any applicable rewriting rules until we are left with an irreducible word---the end result of this rewriting process does not depend on the order in which we applied rewriting rules along the way.

\defi
A rewriting system $(\Sigma, T)$ is called
\enum
\item \emph{length-reducing} if $|R| < |L|$ for every $(L, R) \in T$;
\item \emph{special} if $R = 1$ for every $(L, R)\in T$;
\item \emph{monadic} if $|R| \leq 1$ for every $(L, R) \in T$; and
\item \emph{2-monadic} if $|L| \leq 2$ for every $(L, R) \in T$ and it is length-reducing.
\eenum
\edefi

A finite length-reducing rewriting system is necessarily terminating.   There is a simple algorithm by which one can determine whether or not such a rewriting system is confluent, and hence convergent (see, for example, \cite[Proposition 2.4]{ConfluentThueSystems}).

We shall be concerned with finite, convergent, monadic rewriting systems. Such a rewriting system is called \emph{normalized} if 
$L$ has length at least two and every proper subword of $L$ is reduced for every $(L, R) \in T$.  The following, which is Theorem 1 in \cite{AvenhausMadlener}, shows that we may assume withut loss of generality that our rewriting systems are normalized.

\lem
\label{lem:NormalizedWLOG}
If $(\Sigma, T)$ is a finite, convergent, monadic rewriting system, then there exists a normalized, finite, convergent, monadic rewriting system $(\Sigma', T')$ such that $(\Sigma, T)$ and $(\Sigma', T')$ present isomorphic monoids.
\elem

\section{Potential Obstructions to being Plain}
\label{sec:Obstruction}

For the remainder of the paper, we suppose that $G$ is a group presented by a finite, convergent, monadic rewriting system $(\Sigma, T)$.  By Lemma \ref{lem:NormalizedWLOG} we may assume without loss of generality that $(\Sigma, T)$ is normalized.  

We now show how fundamental results from the 1970's and 1980's combine to allow us to conclude that $G$ may be constructed as the fundamental group of a graph of groups.  Combining important results of Muller and Schupp \cite{MullerSchupp} with those of Dunwoody \cite{Dunwoody} yields that the finitely-generated virtually-free groups are exactly the groups for which the word problem is a context-free language.  Using this characterization, Diekert \cite[Theorem 5]{Diekert} showed that the groups which admit a presentation by a finite, convergent, length-reducing rewriting system form a proper subclass of the virtually-free groups.  Karrass, Pietrowski and Solitar \cite{KPS} characterized the finitely-generated virtually-free groups as the fundamental groups of finite graphs of groups in which the vertex groups are finite.  Thus we have that there exists a finite graph of groups $\Delta$ in which vertex groups are finite and such that $\Delta$ encodes a way to construct a group $\pi(\Delta)$ isomorphic to $G$.

We now interpret the conclusion of the previous paragraph in more detail.  More specifically:
\enum
\item $\Delta$ is a finite, connected, undirected graph with no multi-edges (note that loops are allowed);
\item each vertex $v_i$ is labeled by a finite group $A_i$; and
\item each edge $e$ is labeled by a (necessarily finite) group $K$ and monomorphisms $\phi_{1}: K \to A_i$ and $\phi_{2}: K \to A_j$ into the groups labeling the vertices incident to $e$ (with two monomorphisms into the same vertex group in the case that the edge is a loop).
\eenum
Let $T_0 \subset T_1 \subset \dotsb \subset T_p$ be a sequence of nested subtrees of $\Delta$ such that $T_0$ is a single vertex $\{v_0\}$, $T_p$ is a spanning tree, and each $T_i$ is obtained from $T_{i-1}$ by adding one more vertex $v_i$ and one more edge $e_i$.  Let $e_{p+1}, \dots, e_q$ be the remaining edges in $\Delta$.  Let $K_i$ be the edge group of $e_i$ with monomorphisms $\phi_{i, 1}$ and $\phi_{i, 2}$.  For each $i$ we let $R_i$ be the set comprising all of the relations expressed in the multiplication table for $A_i$ so that $\langle A_i \mid R_i \rangle$ is a finite presentation of $A_i$.  Finally, let $t_{p+1}, \dotsc, t_q$ be new symbols.  The group $\pi(\Delta)$ has a finite presentation $\langle X \mid R\rangle$ with
\[
X = A_0 \cup \dots \cup A_p \cup \{t_{p+1}, \dotsc, t_q\}
\]
and
\begin{align*}
R &= R_0 \cup \dotsb \cup R_p \\ 
&\quad \cup \{\phi_{i, 1}(k) = \phi_{i,2}(k) \text{ for every } 1 \leq i \leq p \text { and every } k \in K_i\} \\
&\quad \cup\{t_i^{-1} \phi_{i, 1} (k) t_i = \phi_{i, 2} (k) \text{ for every } p+1 \leq i \leq q \text { and every } k \in K_i\}.   
\end{align*}
It is important to note that the choices made (for example, the choice of spanning subtree $T_p$) do not affect the isomorphism type of $\pi(\Delta)$.

Without loss of generality we may assume that, for edges that are not loops, the edge homomorphisms $\phi_{i,1}$ and $\phi_{i,2}$ are not surjective (that is, the order of an edge group is strictly less than the order of each vertex group to which the edge is incident). If this were not the case, then we could identify the incident vertices and omit the edge to obtain a more simple graph of groups which presents an isomorphic group.

To prove that $G$ is a plain group, it suffices to show that the edge groups in $\Delta$ are trivial, for in this case the relations associated to edge homomorphisms serve only to identify all of the identity elements from vertex groups, and $\pi(\Delta)$ is isomorphic to the free product of the finite groups $A_0, \dotsc, A_p$ and the free group of rank $q-p$.
To this end we observe some consequences of $\Delta$ having a nontrivial edge group.  We note that $\pi(\Delta)$ may be constructed iteratively using a sequence of free products with amalgamation (one amalgam for each of the edges $e_0, \dotsc, e_p$ in a spanning subtree $T_p$ of $\Delta$) followed by a sequence of HNN extensions (one HNN extension for each of the edges $e_{p+1}, \dotsc, e_q$ not in the spanning subtree).   The following lemma follows from the classical embedding theorems associated to each construction, and the observation that any edge that is not a loop is contained in some spanning subtree of $\Delta$.

\lem
If adjacent vertices are labeled $A_i$ and $A_j$, and the edge group is labeled by $K$, then $G$ contains a subgroup isomorphic to the free product of $A_i$ and $A_j$ with amalgamation over subgroups isomorphic to $K$.  If a vertex is labeled $A$, and a loop at $A$ is labeled $K$, then $G$ contains a subgroup isomorphic to an HNN extension of $A$ with associated subgroups isomorphic to $K$. 
\elem

Our plan is simply to show that $G$ may not contain subgroups of the types described in the lemma.   A nontrivial edge group must be of one of the following types (in each case, we call the edge group $K$):
\enum
\item a loop with cyclic vertex group $A$;
\item an edge with cyclic incident vertex groups $A_i$ and $A_j$;
\item a loop with noncyclic vertex group $A$ such that $|K| = |A|$;
\item a loop with noncyclic vertex group $A$ such that $|K| < |A|$; or
\item an edge with incident vertex groups $A_i$ and $A_j$ which are not both cyclic.
\eenum

The following is a special case of a more general result proved by Madlener and Otto.

\lem
\cite[Theorem 2.3]{MadlenerOtto}
\label{lem:CentralizersOfInfiniteOrderElements}
If $g \in G$ is an element of infinite order, then the centralizer of $g$ in G is isomorphic to $\Z$.
\elem

Madlener and Otto's result can be used to exclude the first three types of nontrivial edge groups.

\lem
\label{lem:BadCentralizer}
The graph of groups $\Delta$ does not contain nontrivial edge groups of type (1), (2), or (3).
\elem

\pf
Suppose that $\Delta$ contains a loop with vertex group $A$, edge group $K$, and homomorphisms $\phi_1, \phi_2 \colon K \to A$.  In case (1), since $A$ is cyclic, $A$ has a unique subgroup of order $|K|$.  In case (3), the maps $\phi_1$ and $\phi_2$ are surjective.  In either of these cases, the images of $\phi_1$ and $\phi_2$ coincide, so $\phi = \phi_2 \circ \phi_1^{-1}$ is an automorphism of $\phi_1(K)$. Now $G$ contains a subgroup isomorphic to 
\[
\langle A, t \mid R, t^{-1} k t = \phi(k) \text{ for all } k \in \phi_1(K) \rangle,
\]
where $R$ comprises the relations expressed in the multiplication table of $A$.  Since $\phi_1(K)$ is a finite group, $\phi^m$ is trivial for some positive integer $m$.  It follows that $t^{-m} k t^m = k$ for all $k \in \phi_1(K)$.  Since $\phi_1(K)$ is nontrivial, this means that the centralizer of $t^m$ (an element of infinite order) contains nontrivial elements of finite order.  This contradicts Lemma \ref{lem:CentralizersOfInfiniteOrderElements}, so $\Delta$ cannot contain a nontrivial edge group of types (1) or (3).

Finally, we consider an edge group with cyclic incident vertex groups $A_i$ and $A_j$, edge group $K$, and homomorphisms $\phi_1 \colon K \to A_i$ and $\phi_2\colon K \to A_j$.  Recall that we assumed without loss of generality that edge groups of non-loops must embed as proper subgroups of the incident vertex groups, so $1 < |K| < \min \{|A_i|, |A_j|\}$. Let $a\in A_i \setminus \phi_1(K)$, $b \in A_j \setminus \phi_2(K)$, and nontrivial $c\in \phi_1(K)$.  Then the infinite order element $ab$ commutes with $c$.  This contradicts Lemma \ref{lem:CentralizersOfInfiniteOrderElements}, so $\Delta$ cannot contain a nontrivial edge group of type (2).
\epf

Nontrivial edge groups of types (4) and (5) are not as easily eliminated, but we can see from the following lemma that the only potential obstruction is an amalgamated product of finite subgroups of $G$:

\lem
\label{lem:EdgeGroupsImplyFPWA}
If $\Delta$ contains a nontrivial edge group of type (4) or (5), then $G$ contains a subgroup isomorphic to a free product with amalgamation
$A \ast_K B$,
where $A$ is a non-cyclic finite group, $B$ is a finite group, and $1 < |K| < \min \{|A|, |B|\}$. 
\elem

\pf
In the case of a type (5) edge group, we clearly do not lose generality by assuming that $A$ is the non-cyclic factor.  In the case of a type (4) edge group with vertex group $A$, edge group $K$, and homomorphisms $\phi_1, \phi_2 \colon K \to A$, we write $\phi = \phi_2 \circ \phi_1^{-1}$ (which, in this case, is an isomorphism from one copy of $K$ in $A$ to another).  Now there exists a subgroup in $G$ presented by
\[
\langle A, t \mid R, t^{-1} k t = \phi(k) \textrm{ for all } k \in \phi_1(k) \rangle
\]
where $R$ comprises the relations expressed in the multiplication tabe of $A$.  The subgroups $t^{-1} A t$ and $A$ generate a subgroup of $G$ which is isomorphic to $(t^{-1} A t) *_K A$.
\epf

To complete our proof of Gilman's conjecture it suffices to show that $G$ cannot contain a free product of finite subgroups, not both cyclic, amalgamated over subgroups which are nontrivial and proper in each factor.  In the next section we show how to identify the finite subgroups of $G$, and we explore their combinatorial and geometric properties.

\section{Finite Order Elements and Subgroups}
\label{sec:FiniteOrder}

We continue to suppose that $G$ is a group presented by a normalized, finite, convergent, monadic rewriting system $(\Sigma, T)$. 

Let $\Gamma$ be the directed Cayley graph of $G$ with respect to $\Sigma$. Thus $\Gamma$ is the labeled directed graph with vertex set $V(\Gamma) = G$, edge set 
\[
E(\Gamma) = \{(g, h) \mid \exists x\in \Sigma, [x] = g^{-1}h \},
\]
and labeling map $L \colon E(\Gamma) \to \Sigma$ defined by $(g, h) \mapsto x$.  We note that, because $(\Sigma, T)$ is normalized, distinct letters represent distinct group elements, and no letter represents the identity.  It follows that $\Gamma$ has no loops or multi-edges.  We also note that rewriting systems of this type have normal forms, that is, for each $U\in \Sigma^*$, there is a unique word $V$ that is shortest among all words in $[U]$, and $V$ is also the unique reduced word in $[U]$.  In $\Gamma$, this means that there are unique geodesic dipaths between any two vertices $g$ and $h$.  For any $g\in G$, we will write $U_g$ for the normal form of $g$ in $\Sigma^*$, and we shall refer to $U_g$ as the reduced representative of $g$.

\begin{figure}[h!]
\begin{tikzpicture}[scale=0.7]

\node (a0) [CayVert,label=-90:{$g=a_0$}] at (0, 0) {};
\node (a1) [CayVert,label=-90:$a_1$] at (2, 0) {};
\node (a2) [CayVert,label=-90:$a_2$] at (4, 0) {};
\node (a5) [CayVert,label=-90:{$a_m=h$}] at (10, 0) {};

\coordinate (a3) at (6, 0);
\coordinate (a4) at (8, 0);

\draw [CayEdge] (a0) to node[midway,above] {$x_1$} (a1);
\draw [CayEdge] (a1) to node[midway,above] {$x_2$} (a2);
\draw [CayEdge] (a2) to node[midway,above] {$x_3$} (a3);
\draw [CayEdge] (a4) to node[midway,above] {$x_m$} (a5);

\node [circle,fill,inner sep=1 pt] at (6.5, 0) {};
\node [circle,fill,inner sep=1 pt] at (7, 0) {};
\node [circle,fill,inner sep=1 pt] at (7.5, 0) {};

\end{tikzpicture}
\caption{The unique geodesic dipath from $g$ to $h$.}
\label{fig:geodesic}
\end{figure}
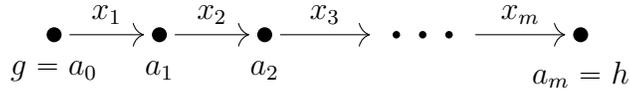

Consider a rewriting rule $(L, R) \in T$.  Let $g$ be a vertex in $\Gamma$, let $\rho_{L}$ be the dipath from $g$ with label $L$, and let $\rho_R$ be the dipath from $g$ with label $R$.  It is clear that $\rho_L$ and $\rho_R$ have the same endpoints---this is simply because the rewriting rules determine equality in the group.   What is characteristic of monadic rewriting systems is the observation that the endpoints of $\rho_L$ are the \emph{only} vertices visited by $\rho_R$.  It follows that if $U, V \in \Sigma^*$ and $U \xrightarrow{*} V$, then the dipath from $g$ with label $V$ visits only vertices visited by the dipath from $g$ with label $U$.  That is:

\lem
\label{lem:PlainGeometry}
Suppose that $g, h \in G$ are distinct and that $x_1 \dotsm x_m \in \Sigma^*$ is the geodesic representative for $g^{-1}h$.  Let $a_0, a_1, \dotsc, a_m$ be the vertices in $\Gamma$ visited by the dipath from $g$ with label $x_1 \dotsm x_m$. (See Figure \ref{fig:geodesic}.)

Then, every dipath from $g$ to $h$ is a concatenation of paths $\rho_1 \rho_2 \dotsm \rho_m$ such that $\rho_i$ is a dipath from $a_{i-1}$ to $a_i$.
\elem
\pf
Consider an arbitrary dipath from $g$ to $h$ corresponding to the word $y_1y_2\dotsm y_n \in \Sigma^*$.  Then $g^{-1}h$ is represented by the word $y_1 y_2 \dotsm y_n$, hence $y_1 y_2 \dots y_n \xrightarrow{*} x_1 x_2 \dotsm x_m$.  Thus there exist words $U_0, \dotsc, U_\ell \in \Sigma^*$ such that
\[
y_1 y_2 \dotsm y_n = U_0 \to U_1 \to \dotsm \to U_\ell = x_1 x_2 \dotsm x_m.
\]
Inductively applying the observation above, $U_\ell$ visits only vertices visited by the path $U_0$.
\epf

An immediate consequence is that once two geodesic dipaths diverge, they cannot rejoin:

\cor
\label{cor:SingleEdge}
If $g \neq h$ and there are two internally disjoint dipaths from $g$ to $h$, then there is a single directed edge from $g$ to $h$.  In particular, suppose that $U$ and $V$ are equivalent words representing $g^{-1}h$.  If $U$ and $V$ begin with different letters and every proper prefix of each is reduced, then  there is a single letter $x$ so that $U \xrightarrow{*} x$ and $V \xrightarrow{*} x$.
\ecor

\pf
The first statement follows immediately from Lemma \ref{lem:PlainGeometry}.  Suppose that $U$ and $V$ satisfy the hypothesis of the second statement.  Let $U'$ be the prefix of $U$ which includes all but the last letter of $U$, and let $V'$ be the prefix of $V$ which includes all but the last letter of $V$.  Since $U'$ and $V'$ are reduced and begin with different first letters, the corresponding paths are internally disjoint and have distinct terminal endpoints.  It follows that $U$ and $V$ are internally disjoint.  
\epf

We use Lemma \ref{lem:PlainGeometry} and its corollary to identify the finite order elements and describe the structure of finite order subgroups of $G$.

\defi
A finite subgroup $A = \{1, a_1, a_2, \dotsc, a_n\}$ of $G$ has the \emph{distinct first letter form} (or \emph{DFL form}) if there exist letters $x_i \in \Sigma$ and words $W_i \in \Sigma^*$ so that the reduced representatives for nontrivial elements of $A$ are
\[
a_i := x_iW_i,
\]
where at least two of the letters $x_i$ are distinct.  (In particular, note that a subgroup in DFL form must have at least two nontrivial elements.) We say that $A$ has the \emph{reduced cyclic form} (or \emph{RC form}) if there exists a word $U \in \Sigma^*$ such that, reordering if necessary, the reduced representatives for the nontrivial elements of $A$ are
\[
a_i := U^i.
\]
\edefi

The next lemma demonstates the profound consequences of the monadic hypothesis.

\lem
\label{lem:TailWord}
Suppose that $A = \{1, a_1, a_2, \dotsc, a_n\}$ has DFL form
\[
a_i := x_iW_i, \quad x_i \in \Sigma, W_i \in \Sigma^*, 1 \leq i \leq n.
\]
Then all of the $x_i$ are distinct, and all of the words $W_i$ are the same word $W$, so that
\[
U_{a_i} = x_iW, \quad 1 \leq i \leq n.
\]
We will refer to the word $W$ as the \emph{tail word for $A$}.
\elem
\pf
Without loss of generality, we may assume that $|W_1|$ is maximal among the lengths $|W_1|, \dotsc, |W_n|$.  Let $W=W_1$.  Since $A$ has DFL form, there is some $i$ so that $x_i \neq x_1$.  Since the dipaths labeled $x_1W$ and $x_iW_i$ emanating from the identity in the Cayley graph are geodesics with distinct first edges, they do not share any vertex other than the identity.

There exists a group element $a_j\in A$ such that $a_i a_j = a_1$.  We have $x_i W_i x_j W_j \xrightarrow{*} x_1W$, and the path $\alpha$ corresponding to $x_1W$ is geodesic, so the path $\beta$ corresponding to $x_i W_i x_j W_j$ must pass through every vertex of the path $\alpha$ by Lemma \ref{lem:PlainGeometry}.  Write $\beta = \rho \sigma$, where $\rho$ is the portion of the path labeled by $x_iW_i$ and $\sigma$ is the portion labeled by $x_j W_j$.

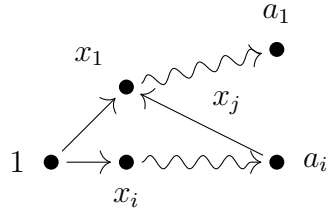
\begin{figure}[h]
\begin{tikzpicture}[scale=0.5]

\node (e) [CayVert,label=180:$1$] at (0,0) {};
\node (x1) [CayVert,label=135:{$x_1$}] at (2, 2) {};
\node (a1) [CayVert,label=90:{$a_1$}] at (6, 3) {};
\node (xi) [CayVert,label=-90:{$x_i$}] at (2, 0) {};
\node (ai) [CayVert,label=0:{$a_i$}] at (6, 0) {};

\draw [CayEdge] (e) to (x1);
\draw [CayEdge] (e) to (xi);
\draw [CayWord] (x1) to (a1);
\draw [CayWord] (xi) to (ai);

\draw [CayEdge] (ai) to node [midway, above right] {$x_j$} (x1);

\end{tikzpicture}
\caption{The path $x_jW_j$ from $a_i$ to $a_1$ must pass through all of the vertices from $x_1$ to $a_1$.  By length arguments, the edge $x_j$ goes directly to the vertex $x_1$, and then $W_j$ follows the path $W$ exactly.}
\label{fig:TailWord}
\end{figure}

The path $\rho$ does not share any vertices other than the vertex 1 with the path $\alpha$, so the rest of the vertices of $\alpha$ must appear along $\sigma$, whose length is $|x_jW_j|$.  It follows from the maximality of $|W|$ that the path $\sigma$ must start with an edge from the vertex $a_i$ to the vertex $x_1$, and it must then follow along $\alpha$ directly to $a_1$.  (If $\sigma$ deviated from $\alpha$ or took longer than a single step to get to $x_1$, then $W_j$ would need to be longer than $W$.)  See Figure \ref{fig:TailWord}.

Now we have $a_1 := x_1 W$ and $a_j := x_j W$. 

For each $\ell$ such that $1 \leq \ell \leq n$, let $S_\ell$ be the set of nontrivial elements in $A$ which have first letter $x_\ell$.  Since $a_1 \neq a_j$, one of the sets $S_1$ or $S_j$ has $n/2$ elements or fewer.  Without loss of generality we may assume that $|S_1| \leq n/2$.  For each $a_k \not \in S_1$, there exists $a_{k'}$ such that $a_k a_{k'} = a_1$. It follows as above that $a_{k'} := x_{k'} W$.   Since $a_{k'} = a_k^{-1} a_1$, each $a_k$ yields a different element $a_{k'}$. We now have that there are at least $n/2 + 1$ elements, including $a_1$, for which the geodesic representative has the form $x_\ell W$.  Without loss of generality we may assume that these elements are $a_1, a_2, \dots, a_m$.  

We note that the sets $S_1, \dots, S_m$ are nonempty and disjoint subsets of $\{a_1, \dots, a_n\}$.  Since $m \geq n/2 + 1$ and $\sum_{i=1}^m |S_i| \leq n$, at least one of the sets has exactly one element.  Without loss of generality we may assume that $|S_1| = 1$.  For each $a_k \neq a_1$, there exists $a_{k'}$ such that $a_k a_{k'} = a_1$.  It follows as above that $a_{k'} := x_{k'} W$.   We have that for each element $g$ of $A$ other than $a_1$, including 1, the geodesic from $g$ to $a_1$ is of the form $x_\ell W$.  Since there are $n$ elements of this form, and $n$ nontrivial elements in $A$, the result is proved.
\epf

\lem
\label{lem:RCorDFL}
Let $A \leq G$ be a finite subgroup.  Then there is $g\in G$, such that one of the following is true:
\enum
\item the conjugate $g^{-1}Ag$ has RC form, or 
\item the conjugate $g^{-1}Ag$ has DFL form.
\eenum
Moreover, if there is some nontrivial element $a\in A$ represented by a word $W = x_1x_2\dotsm x_r$ which has shortest length in the conjugacy class of $a$, then there is some $k \geq 0$ and some (possibly empty) prefix $P$ of $W$ such that $g$ is represented by $W^k P$.
\elem
\pf
The order 2 case is trivial, so assume that $|A| > 2$.  Suppose that the reduced words representing all nontrivial elements begin with the same letter.

Let $a\in A$ be represented by $W = x_1 x_2 \dotsm x_r$. The conjugation $x_1^{-1} A x_1$ cyclically permutes the first letters of each word to the end.  After reducing the words if necessary, if two words have distinct first letters, then $x_1^{-1} A x_1$ has DFL form. Otherwise, repeat the process, rotating the first letter of each word to the end, reducing, and checking for distinct first letters.

After some finite number of steps, we will eventually reach a conjugate $g^{-1} Ag$ of $A$ either with DFL form, or in which all words are cyclically reduced and do not have distinct first letters when continuing to cyclically permute their letters.  In the latter case, it follows that there is some word $U$ so that each nontrivial reduced word of $g^{-1} Ag$ is $U^{e_i}$ for some exponent $e_i$. 

Now suppose that the word $W$ has minimal length in the conjugacy class of $a$.  Then the cyclic conjugates of $W$ are all reduced, so no rewriting rules are applied to the conjugates of $a$ during this procedure.  (Note that the procedure may cyclically conjugate the words of $A$ more than $r$ times, in which case the conjugating letters will repeat.) Therefore, $g$ is represented by a word of the form $W^kP$, where $k \geq 0$ and $P$ is a (possibly empty) prefix of $W$.
\epf

\cor
\label{cor:AllShortest}
Suppose that $A = \{1, a_1, a_2, \dotsc, a_n\}$ is a non-cyclic finite subgroup of $G$.  For each $i$, let $\ell_i$ be the length of the shortest representative among the conjugates of $a_i$.  Then there is a value $\ell$ such that $\ell_i = \ell$ for all $i$.  Moreover, there is a conjugate of $A$ with DFL form whose nontrivial reduced representatives all have length $\ell$.
\ecor
\pf
Without loss of generality, suppose that $\ell_1 = \min_i \{\ell_i\}$, and, replacing $A$ with a conjugate if necessary, suppose that $a_1$ has reduced representative $W$ of length $\ell_1$.  If we apply the cycle-and-reduce procedure outlined in the previous proof, the word $W$ will be cycled but never reduced (by the definition of $\ell_1$).  Since $A$ (and therefore any conjugate of $A$) is not cyclic, this procedure must end with a conjugate $A'$ in DFL form. By Lemma \ref{lem:TailWord}, the nontrivial reduced representatives of $A'$ all have length $\ell_1$.  By the minimality of $\ell_1$, this shows that all of the $\ell_i$ are equal.
\epf

In addition to sharing a tail word, the important properties of subgroups in DFL form are given in the following proposition.

\prop
\label{prop:DFLProperties}
Suppose that $A = \{1, a_1, a_2, \dotsc, a_n\}$ has DFL form, with
\[
a_i := x_i W, \quad 1 \leq i \leq n.
\] 
Then:
\enum
\item Let $V \in \Sigma^*$ such that $[W]^{-1} = [V]$.  If $|W| \geq 1$, then $V\in \Sigma$, and if $W$ is the empty word, then so is $V$.
\eenum
Suppose further that there is some index $j$ such that $x_jW$ has minimal length in its conjugacy class.  Then:
\enum
\setcounter{enumi}{1}
\item For any $i, k$ with $a_ia_j = a_k$, there is a rule $(x_iWx_j, x_k) \in T$.
\item For any $i$ with $a_ia_j = 1$, there is a rule $(x_iWx_j, V) \in T$.
\item The word $Wx_jW$ is reduced.
\eenum
\eprop
\pf
Properties (1) through (4) are trivial if $W$ is the empty word.  Suppose that $W$ is not the empty word.  

\begin{figure}[h]
\begin{tikzpicture}[scale=0.7]

\node (e) [CayVert,label=180:{$1$}] at (0,0) {};
\node (x1) [CayVert,label=135:{$x_1$}] at (2, 2) {};
\node (x2) [CayVert,label=225:{$x_2$}] at (2, -2) {};
\node (a1) [CayVert,label=45:{$a_1$}] at (6, 2) {};
\node (a2) [CayVert,label=-45:{$a_2$}] at (6, -2) {};
\node (g) [CayVert, label=0:$g$] at (4, 0) {};

\draw [CayEdge] (e) to (x1);
\draw [CayEdge] (e) to (x2);
\draw [CayEdge] (a1) to (g);
\draw [CayEdge] (a2) to (g);
\draw [CayWord] (x1) to (a1);
\draw [CayWord] (x2) to (a2);
\draw [CayWord] (g) to (e);

\end{tikzpicture}
\caption{The two dipaths from $1$ to $g$ are internally disjoint, so there must be a single edge from $1$ to $g$.}
\label{fig:Winv}
\end{figure}
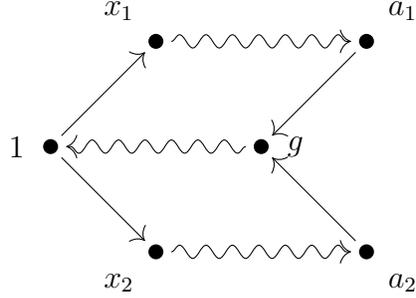

Since $A$ is closed under inverses, there are some $i$ and $j$ so that $a_i = a_1^{-1}$ and $a_j = a_2^{-1}$.  There is a unique path labeled $W$ ending at the vertex $1$, so the paths $x_1 W x_i$ and $x_2 W x_j$ are internally disjoint paths ending at the same vertex $g$ (see Figure \ref{fig:Winv}).  It follows from Corollary \ref{cor:SingleEdge} that there is a single directed edge from $1$ to $g$.  This establishes property (1).  In particular, note that $WV \xrightarrow{*} 1$.

For the remainder of the proof, suppose that there is some index $j$ such that $x_jW$ has minimal length in its conjugacy class.

Let $a_ia_j = a_k$.  Then $x_iWx_jW \xrightarrow{*} x_kW$.  It follows that $x_iWx_jWV \xrightarrow{*} x_kWV$.  Applying the reduction $WV \xrightarrow{*} 1$ to both sides, we have $x_iWx_j \xrightarrow{*} x_k$.  Any proper subword of $x_iWx_j$ is either a subword of $x_iW$ or $Wx_j$.  Both of these words are reduced---the former by assumption, and the latter because it is a conjugate of $x_jW$ which is assumed to have minimal length in its conjugacy class.  Since every proper subword of $x_iWx_j$ is reduced, the reduction $x_iWx_j \xrightarrow{*} x_k$ must be the application of a single rewriting rule $(x_iWx_j, x_k)$.  This establishes property (2).

Similarly, suppose that $a_ia_j = 1$.  Then an analogous argument shows that the reduction $x_iWx_j \xrightarrow{*} V$ must be the application of a single rewriting rule $(x_iWx_j, V)$.  This establishes property (3).

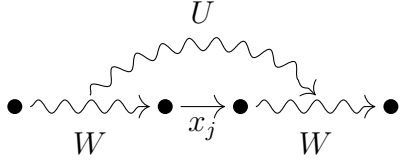
\begin{figure}[h!]
\begin{tikzpicture}[scale=0.5]

\node (x1)   [CayVert] at (0, 0) {};
\node (a1)   [CayVert] at (4, 0) {};
\node (a1x1) [CayVert] at (6, 0) {};
\node (a1a1) [CayVert] at (10, 0) {};

\draw [CayWord] (x1)   to node [midway, below, yshift=-2mm] {$W$}   (a1);
\draw [CayEdge] (a1)   to node [midway, below] {$x_j$} (a1x1);
\draw [CayWord] (a1x1) to node [midway, below, yshift=-2mm] {$W$}   (a1a1);

\draw [CayWord,bend left=45] (2, 0.3) to node[midway, above, yshift=2mm] {$U$} (8, 0.3);

\end{tikzpicture}
\caption{Note that the shortcut path $U$ must connect a vertex from the first $W$ path to a vertex in the second $W$ path.  The two middle vertices cannot be on $U$.}
\label{fig:WxWshortcut}
\end{figure}

Assume that $Wx_jW$ is not reduced.  That is, by Lemma \ref{lem:PlainGeometry}, a path labeled $Wx_jW$ in the Cayley graph has a shortcut (see Figure \ref{fig:WxWshortcut}).

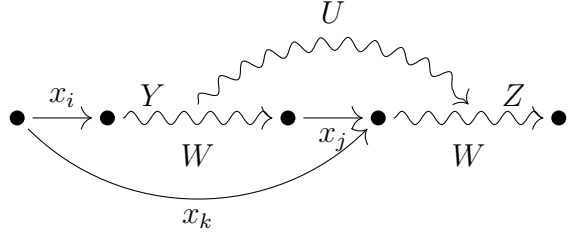
\begin{figure}[h]
\begin{tikzpicture}[scale=0.6]

\node (e)    [CayVert] at (-2, 0) {};
\node (x1)   [CayVert] at (0, 0) {};
\node (a1)   [CayVert] at (4, 0) {};
\node (a1x1) [CayVert] at (6, 0) {};
\node (a1a1) [CayVert] at (10, 0) {};

\draw [CayEdge] (e)    to node [midway, above] {$x_i$} (x1);
\draw [CayWord] (x1)   to node [midway, below, yshift=-2mm] {$W$}   (a1);
\draw [CayEdge] (a1)   to node [midway, below] {$x_j$} (a1x1);
\draw [CayWord] (a1x1) to node [midway, below, yshift=-2mm] {$W$}   (a1a1);

\draw [CayEdge,bend right=45] (e)     to node [midway, below] {$x_k$} (a1x1);
\draw [CayWord,bend left=45] (2, 0.3) to node[midway, above, yshift=2mm] {$U$} (8, 0.3);

\node at (1, 0.5) {$Y$};
\node at (9, 0.5) {$Z$};

\end{tikzpicture}
\caption{The path $x_kW$ is geodesic, so the path $x_iYUZ$ must pass through every vertex along $x_kW$.}
\label{fig:contradiction}
\end{figure}

Choose $i$ such that $a_ia_j \neq 1$.  Then there is a rewriting rule of the form $(x_i W x_j, x_k)$ by (2), thus the subwords $W x_j$ and $x_j W$ are both reduced.  It follows that the path $U$ in Figure \ref{fig:WxWshortcut} does not pass through the middle two vertices.  Let $Y$ represent the initial part of the first $W$ path before $U$, and let $Z$ represent the final part of the second $W$ path after $U$.  Now consider the path $x_iWx_jW$ emanating from 1.  We have the picture shown in Figure \ref{fig:contradiction}.

The path $x_kW$ is a geodesic, so by Lemma \ref{lem:PlainGeometry} the path $x_iYUZ$ passes through every vertex of the path $x_kW$.  But the vertex $x_k$ cannot be on this path, a contradiction.  Therefore $Wx_jW$ must be reduced.
\epf

We observe that Cochet's result \cite{Cochet} follows from this proposition:

\thm
\label{thm:Cochet}
A group $G$ presented by a finite, special, confluent rewriting system $(\Sigma, T)$ is a free product of cyclic groups.
\ethm
\pf
The rewriting system $(\Sigma, T)$ is finite, monadic, and confluent, so $G$ is virtually free by Diekert's result \cite{Diekert}.  If $A$ is any finite subgroup which is not cyclic, replacing $A$ with a conjugate if necessary, we may assume that $A$ has DFL form, and that the reduced representatives of nontrivial elements of $A$ have shortest length in their conjugacy classes.

Since $A$ is not cyclic, $|A| \geq 4$ and $A$ has at least two nontrivial elements which are not inverses.  Given nontrivial elements $a, b\in A$ with $ab \neq 1$, write $a := xW$, $b := yW$, and $ab := zW$.  Proposition \ref{prop:DFLProperties} states that $(xWy, z)$ is a rewriting rule of $T$.  But this contradicts the assumption that $(\Sigma, T)$ is special, so $A$ must be cyclic.

In the language of Section \ref{sec:Obstruction}, $G$ is isomorphic to $\pi(\Delta)$, where $\Delta$ is a graph of groups whose vertex groups are all cylic.  If there are any nontrivial edge groups, they must be of type (1) or type (2), but Lemma \ref{lem:BadCentralizer} says these types of edges cannot occur in $\Delta$.  Thus $\Delta$ has trivial edge groups and cyclic vertex groups, hence $G$ is the free product of cyclic groups.
\epf

Finally, we explore some consequences for elements in finite subgroups of DFL or RC form whose reduced representatives are minimal length in their conjugacy class.

\lem
\label{lem:NotRC}
Suppose that $A$ is a finite cyclic subgroup of $G$ having order $m+1$ at least 3.  Suppose that there is some $z\in \Sigma$, $W\in \Sigma^*$, and $c\in A$ of order 2 such that $c := zW$ and the word $zWz$ is not reduced.  Then $A$ cannot have RC form.
\elem
\pf
Suppose that $A$ has RC form, generated by an element $g$ of order $m+1$.  Let $U$ be the reduced representative of $g$, so that the reduced representatives for $A$ are:
\[
A = \{1, U, U^2, \dotsc, U^m\}.
\]
Since $c$ has order 2, $m+1$ must be even and $c = g^{(m+1)/2}$.  Then $zW = U^{(m+1)/2}$.  In particular, $U$ begins with the letter $z$, so we can write $U = zV$.  Consider the element $cg = g^{(m+1)/2 + 1}$.  Since $m+1 > 2$, we have
\[
\frac{m+1}{2} + 1 < m + 1,
\]
so $cg \neq 1$.  

According to the RC form for $A$, the reduced form for $cg$ should be:
\[
U^{(m+1)/2 + 1} = U^{(m+1)/2}U = (zW)(zV).
\]
But this is not reduced, since it contains $zWz$ as a subword.  Therefore $A$ cannot have RC form.
\epf

For the sake of clarity, we introduce the following terminology:

\defi
The word $xW$ is \emph{appended first letter reducible} (or \emph{AFL-reducible}) if the word $xWx$ is reducible.  The word $xW$ is \emph{appended first letter irreducible} (or \emph{AFL-irreducible}) if the word $xWx$ is irreducible.
\edefi

\lem
\label{lem:UniqueCyclicConjugate}
Suppose that $A$ is a finite subgroup of $G$ having DFL form.  Let $a\in A$ be nontrivial and have reduced representative $xW$ for some $x\in \Sigma$, $W\in \Sigma^*$.  Suppose that $xW$ has shortest length among representatives of conjugates of $a$.  Then $xW$ is AFL-reducible, and every other cyclic conjugate of $xW$ is AFL-irreducible.
\elem
\pf
Proposition \ref{prop:DFLProperties} part (2) or (3) shows that $xW$ is AFL-reducible.  On the other hand, by part (4) of that proposition, the word $WxW$ is reduced.  Any other cyclic conjugate of $xW$ followed by its first letter is a subword of $WxW$, hence it is reduced.
\epf

\cor
\label{cor:ConjugateNotDFL}
Suppose that $A$ is a finite subgroup of $G$ having DFL form.  Suppose that some nontrivial element $a\in A$ has reduced representative $xW$ which is minimal length in its conjugacy class.  Let $B$ be any finite subgroup containing $a$. Then $g^{-1}Bg$ cannot have DFL form if $g$ is represented by a nontrivial and proper prefix of $xW$.
\ecor
\pf
The word $xW$ is AFL-reducible, and every other cyclic conjugate of $xW$ is AFL-irreducible. Suppose that $g^{-1}Bg$ has DFL form and $g$ is a nontrivial and proper prefix of $xW$.  A reduced representative of $g^{-1}ag$ is a cyclic conjugate of $xW$, and it would have to be AFL-reducible by Proposition \ref{prop:DFLProperties} part (2) or (3), a contradiction.
\epf

\section{Main Result}
\label{sec:Result}

Suppose that $G$ is a group presented by a finite, convergent, monadic rewriting system $(\Sigma, T)$.  By Lemma \ref{lem:NormalizedWLOG} we may assume without loss of generality that $(\Sigma, T)$ is normalized.  

As laid out in Section \ref{sec:Obstruction}, to complete our proof of Gilman's conjecture it remains to show that $G$ cannot contain a subgroup isomorphic to $A *_C B$, where $A$ and $B$ are finite subgroups, $A$ is non-cyclic, and $C$ is a nontrivial, proper subgroup of $A$ and of $B$.  We shall proceed by showing that $A$ and $B$ cannot both have DFL form, and then we prove that we may replace $A *_C B$ with a conjugate in which $A$ and $B$ are both in DFL form.

\thm
\label{thm:AmalgamatedDFLs}
The group $G$ does not contain a subgroup isomorphic to a group
\[
A \ast_C B
\]
where $A$ and $B$ both have DFL form, and $1 < |C| < \min \{|A|, |B|\}$.
\ethm
\pf
Suppose that $A *_C B$ is such an amalgamated product.  Let $a \in A \setminus C$, $b \in B \setminus C$ and $c \in C \setminus \{1\}$.  By the normal form for free products with amalgamation, $ab$ has infinite order.  We shall show that the elements in 
\[
\{(ab)^k \mid k \in \Z\}
\]
are represented by reduced words of uniformly bounded length.  Since there are only finitely many such words, this contradicts $ab$ having infinite order, establishing the theorem.

Let $A = \{1, a_1, a_2, \dotsc, a_n\}$ and $B = \{1, b_1, b_2, \dotsc, b_m\}$.  Write the reduced representatives of $A$ and $B$ as $a_i := x_iW$ and $b_i := y_iW$ respectively.  (Note that the tail words must be the same, since $A$ and $B$ share nontrivial elements from $C$.)  Without loss of generality, suppose that $c = a_1 = b_1$, $a = a_2$, and $b = b_2$.  Let $z = x_1 = y_1$, so that $c$ is represented by the word $zW$.

Since $A$ is a group, there is some $a_i$ such that $a_ic = a$.  Since $c$ is an element of $C$ but $a$ is not, $a_i \neq c$.  Since $A$ is a group and $c \neq 1$, $a_i \neq a$.  Without loss of generality (reindexing if necessary), we may suppose that $i = 3$.  Similarly, we may suppose $cb = b_3$. It follows that $x_3 W z \xrightarrow{*} x_2$ and $z W y_2 \xrightarrow{*} y_3$.  Now we have $x_2 W y_2 \xleftarrow{*} x_3 W z W y_2 \xrightarrow{*} x_3 W y_3$.

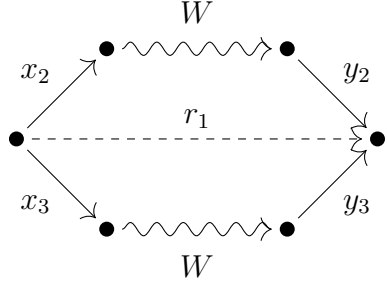
\begin{figure}[h]
\begin{tikzpicture}[scale=0.6]

\node (e)  [CayVert] at (0,0) {};
\node (x2) [CayVert] at (2,2) {};
\node (x3) [CayVert] at (2,-2) {};
\node (a)  [CayVert] at (6,2) {};
\node (a3) [CayVert] at (6,-2) {};
\node (r)  [CayVert] at (8,0) {};

\draw [CayEdge] (e) to node[midway,above left] {$x_2$} (x2);
\draw [CayEdge] (e) to node[midway,below left] {$x_3$} (x3);
\draw [CayWord] (x2) to node[midway,above,yshift=2mm] {$W$} (a);
\draw [CayWord] (x3) to node[midway,below,yshift=-2mm] {$W$} (a3);
\draw [CayEdge] (a) to node[midway,above right] {$y_2$} (r);
\draw [CayEdge] (a3) to node[midway,below right] {$y_3$} (r);

\draw [CayDashEdge] (e) to node[midway,above] {$r_1$} (r);

\end{tikzpicture}
\caption{The top and bottom paths from left to right are internally disjoint, since both $x_2W$ and $x_3W$ are reduced words.}
\label{fig:amalgamatedr1}
\end{figure}

The paths labeled $x_2Wy_2$ and $x_3Wy_3$ are internally disjoint (because $x_2W$ and $x_3W$ are geodesics with $x_2 \neq x_3$) and they are not loops (because $y_2$ is not the inverse of $x_2 W$), so by Lemma \ref{lem:PlainGeometry} there is some $r_1\in \Sigma$ so that $x_2Wy_2 \xrightarrow{*} r_1$ and $x_3Wy_3 \xrightarrow{*} r_1$.  (See Figure \ref{fig:amalgamatedr1}.)  By a similar argument, there is some $s_1 \in \Sigma$ such that $y_2 W x_2 \xrightarrow{*} s_1$.

\begin{figure}[h]
\begin{tikzpicture}[scale=0.6]

\node (e)  [CayVert] at (0,0) {};
\node (x2) [CayVert] at (2,2) {};
\node (a)  [CayVert] at (6,2) {};
\node (r)  [CayVert,label=0:$g$] at (8,0) {};
\node (v)  [CayVert] at (4,0) {};
\node (ab)  [CayVert] at (4,-4) {};
\node (U)  [CayVert] at (6,-2) {};

\draw [CayEdge] (e) to node[midway,above left] {$x_2$} (x2);
\draw [CayWord] (x2) to node[midway,above,yshift=2mm] {$W$} (a);
\draw [CayEdge] (a) to node[midway,above right] {$y_2$} (r);
\draw [CayEdge] (a) to node[midway,below right] {$x_i$} (v);
\draw [CayWord] (v) to node[pos=0.4,above] {$W$} (e);
\draw [CayWord] (r) to node[pos=0.2,below right] {$U$} node[pos=0.8,below right] {$V$} (ab);
\draw [CayEdge] (e) to node[midway,below] {$t$} (U);

\end{tikzpicture}
\caption{The path $y_2U$ is geodesic, so the path $x_iWt$ must pass through vertex $g$, but this is impossible.}
\label{fig:rWreduced}
\end{figure}
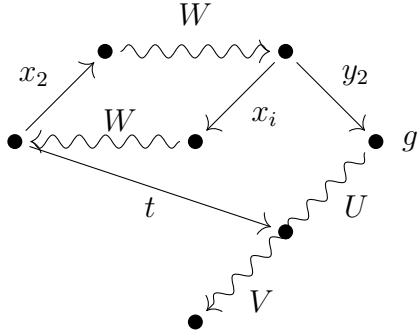

We claim that the word $r_1W$ is reduced.  If not, then $r_1W \to V$ or $r_1W \to tV$, where $W = UV$.  Consider first the case that $r_1 W \to V$.  
It follows that $x_2 W y_2 U \xrightarrow{\ast} 1$.  Then $y_2 U$ spells the inverse of $x_2 W$.  Since  $|y_2 U| < |x_2 W|$, this contradicts the fact that  every nontrivial element in a DFL group has the same length.
Now consider the case that $r_1 W \to tV$.  Let $a_i = a^{-1}$.  We note that $x_i \neq y_2$, since $a_i$ is not in $B$.  Now consider the paths labeled $y_2U$ and $x_iWt$ in Figure \ref{fig:rWreduced}.  By Lemma \ref{lem:PlainGeometry}, the path $x_iWt$ must pass through every vertex of the path $y_2U$, since the latter path is geodesic.  However, $x_iWt$ cannot pass through the vertex labeled $g$ in the figure, otherwise the edge $y_2$ would provide a shortcut on a path that is supposed to be geodesic.  This is a contradiction, so the word $r_1W$ must be reduced.  By a similar argument, the word $s_1W$ is reduced.

\begin{figure}[h]
\begin{tikzpicture}[scale=0.6]

\node (e)  [CayVert] at (0,0) {};
\node (x2) [CayVert] at (2,2) {};
\node (a)  [CayVert] at (6,2) {};
\node (v)  [CayVert] at (4,0) {};

\draw [CayEdge] (e) to node[midway,above left] {$x_2$} (x2);
\draw [CayWord] (x2) to node[midway,above,yshift=2mm] {$W$} (a);
\draw [CayEdge] (a) to node[midway,below right] {$x_i$} (v);
\draw [CayWord] (v) to node[pos=0.4,above] {$W$} (e);
\draw [CayEdge] (x2) to node[pos=0.7,above right] {$r_1$} (v);

\end{tikzpicture}
\caption{The word $Wx_i$ is reduced, so there cannot be a shortcut $r_1$ as in the figure.}
\label{fig:r1Wx2}
\end{figure}
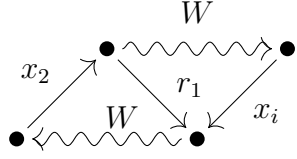

Now consider the word $x_2 W y_2 W x_2$.  We have 
\[
r_1 W x_2 \xleftarrow{*} x_2 W y_2 W x_2 \xrightarrow{*} x_2 W s_1.
\]
We observe that $r_1W$ cannot be the inverse of $x_2$: otherwise, any path labeled $r_1Wx_2$ forms a loop, so there would have to be an edge labeled $r_1$ as in Figure \ref{fig:r1Wx2}.  If $W$ is the empty word, this figure would show $r_1 = x_i$, which is impossible (since then $ab \in A$).  If $W$ is not the empty word, then $Wx_i$ is reduced, so the path $r_1$ cannot provide a shortcut.

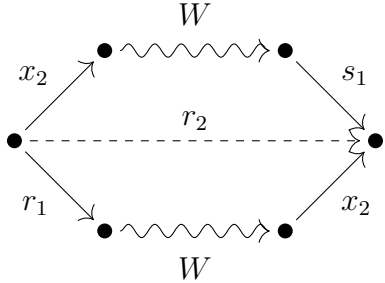
\begin{figure}[h]
\begin{tikzpicture}[scale=0.6]

\node (e)  [CayVert] at (0,0) {};
\node (x2) [CayVert] at (2,2) {};
\node (x3) [CayVert] at (2,-2) {};
\node (a)  [CayVert] at (6,2) {};
\node (a3) [CayVert] at (6,-2) {};
\node (r)  [CayVert] at (8,0) {};

\draw [CayEdge] (e) to node[midway,above left] {$x_2$} (x2);
\draw [CayEdge] (e) to node[midway,below left] {$r_1$} (x3);
\draw [CayWord] (x2) to node[midway,above,yshift=2mm] {$W$} (a);
\draw [CayWord] (x3) to node[midway,below,yshift=-2mm] {$W$} (a3);
\draw [CayEdge] (a) to node[midway,above right] {$s_1$} (r);
\draw [CayEdge] (a3) to node[midway,below right] {$x_2$} (r);

\draw [CayDashEdge] (e) to node[midway,above] {$r_2$} (r);

\end{tikzpicture}
\caption{The top and bottom paths from left to right are internally disjoint, since both $x_2W$ and $r_1W$ are reduced words.}
\label{fig:amalgamatedr2}
\end{figure}

We also have that $r_1W$ and $x_2W$ are both reduced, and they start with distinct letters (otherwise $r_1W$ represents the group element $ab$, and this cannot be equal to $a$).  Therefore the paths across the top and bottom of Figure \ref{fig:amalgamatedr2} are internally disjoint and they are not closed. By Lemma \ref{lem:PlainGeometry}, there is some $r_2 \in \Sigma$ such that $x_2 W y_2 W x_2 \xrightarrow{*} r_2$, and the word $r_2W$ is reduced by the same argument as above.  Similarly, there is some $s_2 \in \Sigma$ such that $y_2 W x_2 W y_2 \xrightarrow{*} s_2$, and $s_2W$ is reduced.

Continuing inductively, we have that the word $(x_2 W y_2 W)^k x_2$ reduces to the single letter $r_{2k}$ for each $k \in \N$. It follows that elements in 
\[
\{(ab)^k \mid k \in \Z\}
\]
are represented by reduced words of bounded length.  But there are only finitely many words of length up to a particular bound, so this implies that $ab$ has finite order, a contradiction.  Therefore $G$ cannot contain a subgroup of the form $A *_C B$.
\epf

\thm
\label{thm:AmalgamatedProducts}
The group $G$ does not contain a subgroup isomorphic to a group
\[
A \ast_C B
\]
where $A$ is a non-cyclic finite group, $B$ is a finite group, and $1 < |C| < \min \{|A|, |B|\}$.
\ethm
\pf
Suppose that $A *_C B$ is such an amalgamated product.  By replacing $A *_C B$ by a conjugate if necessary, we may assume that $A$ has DFL form (Lemma \ref{lem:RCorDFL}) and the nontrivial elements of $A$ are each represented by words which are shortest among all representatives of conjugates of nontrivial elements in $A$ (Corollary \ref{cor:AllShortest}).  We will write these representatives as
\[
A = \{1, zW, x_2W, x_3W, \dotsc, x_nW\}, 
\]
where $z, x_2, x_3, \dotsc, x_n \in \Sigma$, $W\in \Sigma^*$, $c := zW$ and $a := x_2W$.  We shall show that, replacing $A *_C B$ by a further conjugate if necessary, $B$ is also in DFL form.  By Theorem \ref{thm:AmalgamatedDFLs}, this is impossible, completing the proof.

In the case that $C$ has order at least three, then it contains two elements which start with different letters, because $A$ has DFL form.  But these elements are also in $B$, so $B$ has DFL form.

For the remainder of the argument, assume that $C$ (and therefore $c$) has order 2.  By Lemma \ref{lem:RCorDFL}, some conjugate $B' = g^{-1}Bg$ has either RC or DFL form.  Moreover, since $zW$ has minimal length in its conjugacy class, there is some power $k \geq 0$ and a (possibly empty) prefix $P$ of $zW$ so that $g$ is represented by $(zW)^kP$.  Since $c$ has order 2, $k$ must be either 0 or 1.  We now consider subcases depending on whether $P$ is the empty word and whether $g^{-1} B g$ has RC or DFL form.

Suppose $P$ is the empty word, so that either $g = 1$ or $g = c$.  If $g = 1$, then $B$ itself has RC or DFL form.  Note that $zWz$ is not reduced by Lemma \ref{prop:DFLProperties} part (3), so Lemma \ref{lem:NotRC} applies and shows $B$ cannot have RC form, so it must have DFL form.  On the other hand, if $g = c$, then we can write $B' = c^{-1}Bc$, and consider the amalgamated product $c^{-1}(A*_C B)c = A *_C B'$.  Thus we have reduced to the case in which $g = 1$, where we have already concluded $B'$ has DFL form.

Now suppose $P$ is not the empty word, so that $g^{-1}cg$ is a nontrivial cyclic conjugate of $zW$. This conjugate is AFL-irreducible by Lemma \ref{lem:UniqueCyclicConjugate}.  If $B'$ had DFL form, then Proposition \ref{prop:DFLProperties} part (2) or (3) would imply that $g^{-1}cg$ is AFL-reducible, a contradiction.
Therefore $B'$ must have RC form.  The property of having RC form is preserved by cyclic conjugation, so we may replace $g^{-1}Bg$ with a further conjugate $B'' = h^{-1}Bh$, where $h$ is a power of $c$.  We thereby reduce to the case that $P$ is the empty word, therefore $B''$ has DFL form.

We have shown that one of the amalgamated products $A *_C B$, $A *_C B'$, or $A *_C B''$ must satisfy the assumptions of Theorem \ref{thm:AmalgamatedDFLs}, therefore $G$ cannot contain $A *_C B$ as a subgroup.
\epf

Following the discussion in Section \ref{sec:Obstruction}, this completes our main result:

\thm \label{thm:GilmansConjecture}
A group $G$ can be presented by a finite, convergent, monadic rewriting system if and only if $G$ is a plain group.
\ethm

\bibliographystyle{alpha}
\bibliography{GilmansConjectureBib}

\newpage

\end{document}